\documentclass[12pt]{amsart}
\usepackage{amsfonts,amssymb,amscd,amsmath,enumerate,verbatim,calc,latexsym,pstcol,pst-plot,pst-3d,}
\def\NZQ{\mathbb}               

\def\QQ{{\NZQ Q}}
\def\ZZ{{\NZQ Z}}
\def\RR{{\NZQ R}}

\def\frk{\mathfrak}               

\def\Phi{{\frk N}}
\def\a{\alpha}

\def\opn#1#2{\def#1{\operatorname{#2}}} 
\opn\chara{char} 
\opn\length{\ell} 
\opn\pd{pd} 
\opn\rk{rk}
\opn\projdim{proj\,dim} 
\opn\injdim{inj\,dim} 
\opn\rank{rank}
\opn\depth{depth} 
\opn\grade{grade} 
\opn\height{height}
\opn\embdim{emb\,dim} 
\opn\codim{codim}

\opn\Tr{Tr} 
\opn\bigrank{big\,rank}
\opn\superheight{superheight}
\opn\lcm{lcm}
\opn\trdeg{tr\,deg}
\opn\reg{reg} 
\opn\lreg{lreg} 
\opn\ini{in} 
\opn\lpd{lpd}
\opn\size{size}
\opn\mult{mult}
\opn\dist{dist}
\opn\cone{cone}
\opn\lex{lex}
\opn\rev{rev}
\opn\div{div} \opn\Div{Div} \opn\cl{cl} \opn\Cl{Cl}
\opn\Spec{Spec} \opn\Supp{Supp} \opn\supp{supp} \opn\Sing{Sing}
\opn\Ass{Ass} \opn\Min{Min}
\opn\Ann{Ann} \opn\Rad{Rad} \opn\Soc{Soc}
\opn\Syz{Syz} \opn\Im{Im} \opn\Ker{Ker} \opn\Coker{Coker}
\opn\Am{Am} \opn\Hom{Hom} \opn\Tor{Tor} \opn\Ext{Ext}
\opn\End{End} \opn\Aut{Aut} \opn\id{id} \opn\ini{in}

\opn\nat{nat}
\opn\pff{pf}
\opn\Pf{Pf} \opn\GL{GL} \opn\SL{SL} \opn\mod{mod} \opn\ord{ord}
\opn\Gin{Gin}
\opn\Hilb{Hilb}\opn\adeg{adeg}\opn\std{std}\opn\ip{infpt}
\opn\Pol{Pol}
\opn\sat{sat}
\opn\Var{Var}
\opn\Gen{Gen}

\opn\aff{aff} \opn\con{conv} \opn\relint{relint} \opn\st{st}
\opn\lk{lk} \opn\cn{cn} \opn\core{core} \opn\vol{vol}
\opn\link{link} \opn\star{star}
\opn\gr{gr}

\def\Hc{{\mathcal H}}

\def\Fc{{\mathcal F}}
\def\Pc{{\mathcal P}}
\def\Xc{{\mathcal Q}}
\def\Oc{{\mathcal O}}
\def\pot#1#2{#1[\kern-0.28ex[#2]\kern-0.28ex]}

\opn\dirlim{\underrightarrow{\lim}}
\opn\inivlim{\underleftarrow{\lim}}
\let\union=\cup

\let\to=\rightarrow

\def\Implies{\ifmmode\Longrightarrow \else
        \unskip${}\Longrightarrow{}$\ignorespaces\fi}
\def\implies{\ifmmode\Rightarrow \else
        \unskip${}\Rightarrow{}$\ignorespaces\fi}
\def\iff{\ifmmode\Longleftrightarrow \else
        \unskip${}\Longleftrightarrow{}$\ignorespaces\fi}

\let\:=\colon
\newtheorem{Theorem}{Theorem}[section]
\newtheorem{Lemma}[Theorem]{Lemma}
\newtheorem{Corollary}[Theorem]{Corollary}

\newtheorem{Remark}[Theorem]{Remark}

\newtheorem{Example}[Theorem]{Example}

\newtheorem{Definition}[Theorem]{Definition}

\let\epsilon\varepsilon
\let\phi=\varphi
\let\kappa=\varkappa
\def\qed{\ifhmode\textqed\fi
      \ifmmode\ifinner\quad\qedsymbol\else\dispqed\fi\fi}
\def\textqed{\unskip\nobreak\penalty50
       \hskip2em\hbox{}\nobreak\hfil\qedsymbol
       \parfillskip=0pt \finalhyphendemerits=0}
\def\dispqed{\rlap{\qquad\qedsymbol}}

\opn\dis{dis}
\opn\height{height}
\opn\dist{dist}
\def\pnt{{\raise0.5mm\hbox{\large\bf.}}}

\opn\Lex{Lex}

\begin{document}
\title{Smooth Fano polytopes arising from finite partially ordered sets}
\author{Takayuki Hibi and Akihiro Higashitani}
\address{Takayuki Hibi,
Department of Pure and Applied Mathematics,
Graduate School of Information Science and Technology,pe
Osaka University,
Toyonaka, Osaka 560-0043, Japan}
\email{hibi@math.sci.osaka-u.ac.jp}
\address{Akihiro Higashitani,
Department of Pure and Applied Mathematics,
Graduate School of Information Science and Technology,
Osaka University,
Toyonaka, Osaka 560-0043, Japan}
\email{sm5037ha@ecs.cmc.osaka-u.ac.jp}
\begin{abstract}
Gorenstein Fano polytopes arising from finite
partially ordered sets will be introduced.
Then we study the problem of which 
partially ordered sets yield
smooth Fano polytopes.   
\end{abstract}
\maketitle

\section*{Introduction}
An integral (or lattice) polytope is a convex polytope
all of whose vertices have integer coordinates. 
Let $\Pc \subset \RR^d$ be an integral convex polytope 
of dimension $d$.
\begin{itemize}
\item
We say that $\Pc$ is a {\em Fano polytope} if
the origin of $\RR^d$ is the unique integer point
belonging to the interior of $\Pc$. 
\item
A Fano polytope $\Pc$ is called {\em terminal} if
each integer point belonging to the boundary of
$\Pc$ is a vertex of $\Pc$.
\item
A Fano polytope is called {\em Gorenstein} 
if its dual polytope is integral.
(Recall that the dual polytope $\Pc^\vee$
of a Fano polytope $\Pc$ is the convex polytope
which consists of those $x \in \RR^d$
such that $\langle x, y \rangle \leq 1$ for all
$y \in \Pc$, where $\langle x, y \rangle$
is the usual inner product of $\RR^d$.)
\item
A {\em $\QQ$-factorial Fano polytope} is a simplicial Fano polytope,
i.e., a Fano polytope each of whose faces is a simplex. 
\item
A {\em smooth Fano polytope} is a Fano polytope such that
the vertices of each facet form a $\ZZ$-basis of $\ZZ^d$. 
(Sometimes, smooth polytopes denote simple polytopes, 
which are dual polytopes of simplicial polytopes.) 
\end{itemize}
Thus in particular a smooth Fano polytope is 
$\QQ$-factorial, Gorenstein and terminal.

{\O}bro \cite{Oeb} succeeded in finding an algorithm 
which yields the classification list of 
the smooth Fano polytopes for given $d$. 
It is proved in Casagrande \cite{Cas} that
the number of vertices of a Gorenstein $\QQ$-factorial
Fano polytope is at most $3d$ if $d$ is even, and
at most $3d-1$ if $d$ is odd.
B. Nill and M. {\O}bro \cite{NO} classified
the Gorenstein $\QQ$-factorial
Fano polytopes of dimension $d$ with $3d-1$ vertices.
Gorenstein Fano polytopes are classified when $d \leq 4$ 
by Kreuzer and Skarke \cite{KS98}, \cite{KS00} 
and mirro symmetry is studied as the relevance of 
Gorenstein Fano polytopes by Batyrev \cite{Bat94}. 
The study on the classification of terminal 
or canonical Fano polytopes 
was done by Kasprzyk \cite{Kas06}, \cite{Kas08}.
The combinatorial conditions for what it implies 
to be terminal and canonical are explained in 
Reid \cite{Rei83}. 

In the present paper, 
given a finite partially ordered set $P$ 
we associate a terminal Fano polytope
$\Xc_P$.  By using the theory of totally unimodular
matrices, it turns out that these Fano polytopes
are Gorenstein. 
Then we study the problem of which 
partially ordered sets yield
$\QQ$-factorial Fano polytopes.
Finally, it turns out that the Fano polytope
$\Xc_P$ is smooth if and only if
$\Xc_P$ is $\QQ$-factorial.

\section{Fano polytopes arising from finite partially ordered sets}
Let $P = \{ y_1, \ldots, y_d \}$ be a finite partially ordered set
and
\[
\hat{P} = P \cup \{ \hat{0},\hat{1} \},
\]
where $\hat{0}$ (resp. $\hat{1}$) is a unique minimal
 (resp. maximal) element of $\hat{P}$
with $\hat{0} \not\in P$ (resp. $\hat{1} \not\in P$). 
Let $y_0 = \hat{0}$ and $y_{d+1} = \hat{1}$.
We say that $e = \{ y_i, y_j \}$, 
where $0 \leq i, \, j \leq d + 1$
with $i \neq j$, is an {\em edge} of $\hat{P}$
if $e$ is an edge of the Hasse diagram of $\hat{P}$.
(The Hasse diagram of a finite partially ordered set 
can be regarded as a finite nondirected graph.) 
In other words, 
$e = \{ y_i, y_j \}$ is an edge of $\hat{P}$
if $y_i$ and $y_j$ are comparable in $\hat{P}$, say, $y_i < y_j$,
and there is no $z \in P$ with $y_i < z < y_j$.

\begin{Definition}
{\em
Let 
$\hat{P} = \{ y_0, y_1, \ldots, y_d, y_{d+1} \}$ 
be a finite partially ordered set with 
$y_0 = \hat{0}$ and $y_{d+1} = \hat{1}$. 
Let ${\bf e}_i$ denote
the $i$th
canonical unit coordinate vector of $\RR^d$.
Given an edge $e=\{y_i, y_j\}$ of $\hat{P}$ 
with $y_i < y_j$, we define $\rho(e) \in \RR^d$ by setting 
\begin{eqnarray*}
\rho(e)=
\begin{cases}
\, \, \, \, \, \, {\bf e}_i
\;\;\;\;\;\;\;\;&
\text{if} 
\; \; \; j=d+1, \\
\, - {\bf e}_j 
&
\text{if} 
\; \; \;  i=0, \\
\, {\bf e}_i - {\bf e}_j 
&
\text{if} 
\; \; \; 1 \leq i, \, j \leq d.
\end{cases}
\end{eqnarray*} 
Moreover, we write $\Xc_P \subset \RR^d$ for the convex hull 
of the finite set 
\[
\{\,\rho(e) : \, \text{$e$ is an edge of} \;\; \hat{P}\,\}. 
\]
}
\end{Definition}

\begin{Example}
{\em
Let $P = \{ y_1, y_2, y_3 \}$ 
be the finite partially ordered set 
with the partial order $y_1 < y_2$.  
Then ${\hat P}$ together with $\rho(e)$'s
and $\Xc_P$ are drawn below: }
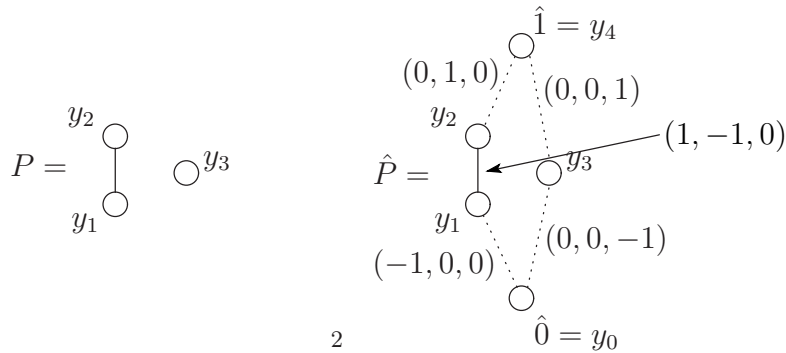
\begin{center}
\unitlength 0.1in
\begin{picture}( 34.1900, 16.3200)(  7.8000,-22.3200)
%
\special{pn 8}%
\special{ar 1702 1488 64 64  0.0000000 6.2831853}%
%
\special{pn 8}%
\special{ar 1330 1296 64 64  0.0000000 6.2831853}%
%
\special{pn 8}%
\special{ar 1330 1652 64 64  0.0000000 6.2831853}%
%
\special{pn 8}%
\special{pa 1330 1358}%
\special{pa 1330 1590}%
\special{fp}%
\put(7.8000,-13.8800){\makebox(0,0)[lt]{$P=$}}%
%
\put(12.4600,-6.0000){\makebox(0,0)[lb]{}}%
\put(10.8000,-11.4400){\makebox(0,0)[lt]{$y_2$}}%
\put(11.0200,-16.9900){\makebox(0,0)[lt]{$y_1$}}%
\put(17.9000,-13.7700){\makebox(0,0)[lt]{$y_3$}}%
%
\special{pn 8}%
\special{ar 3600 1488 64 64  0.0000000 6.2831853}%
%
\special{pn 8}%
\special{ar 3228 1296 64 64  0.0000000 6.2831853}%
%
\special{pn 8}%
\special{ar 3228 1652 64 64  0.0000000 6.2831853}%
%
\special{pn 8}%
\special{pa 3228 1358}%
\special{pa 3228 1590}%
\special{fp}%
\put(26.7800,-13.8800){\makebox(0,0)[lt]{$\hat{P}=$}}%
\put(29.7800,-11.4400){\makebox(0,0)[lt]{$y_2$}}%
\put(30.0000,-16.9900){\makebox(0,0)[lt]{$y_1$}}%
\put(36.8800,-13.7700){\makebox(0,0)[lt]{$y_3$}}%
%
\special{pn 8}%
\special{ar 3456 822 64 64  0.0000000 6.2831853}%
%
\special{pn 8}%
\special{ar 3456 2144 64 64  0.0000000 6.2831853}%
\put(35.1100,-6.1100){\makebox(0,0)[lt]{$\hat{1}=y_4$}}%
\put(35.2200,-22.3200){\makebox(0,0)[lt]{$\hat{0}=y_0$}}%
%
\special{pn 8}%
\special{pa 3422 878}%
\special{pa 3266 1234}%
\special{dt 0.045}%
\special{pa 3500 878}%
\special{pa 3600 1400}%
\special{dt 0.045}%
\special{pa 3256 1710}%
\special{pa 3434 2088}%
\special{dt 0.045}%
\special{pa 3488 2088}%
\special{pa 3600 1566}%
\special{dt 0.045}%
\put(33.4400,-8.8900){\makebox(0,0)[rt]{$(0,1,0)$}}%
\put(35.6600,-9.7700){\makebox(0,0)[lt]{$(0,0,1)$}}%
\put(35.7700,-17.5400){\makebox(0,0)[lt]{$(0,0,-1)$}}%
\put(33.2200,-18.8800){\makebox(0,0)[rt]{$(-1,0,0)$}}%
%
\special{pn 8}%
\special{pa 4178 1288}%
\special{pa 3278 1478}%
\special{fp}%
\special{sh 1}%
\special{pa 3278 1478}%
\special{pa 3348 1484}%
\special{pa 3330 1466}%
\special{pa 3340 1444}%
\special{pa 3278 1478}%
\special{fp}%
\put(41.9900,-12.1100){\makebox(0,0)[lt]{$(1,-1,0)$}}%
\end{picture}%
\end{center}
\begin{center}
\unitlength 0.1in
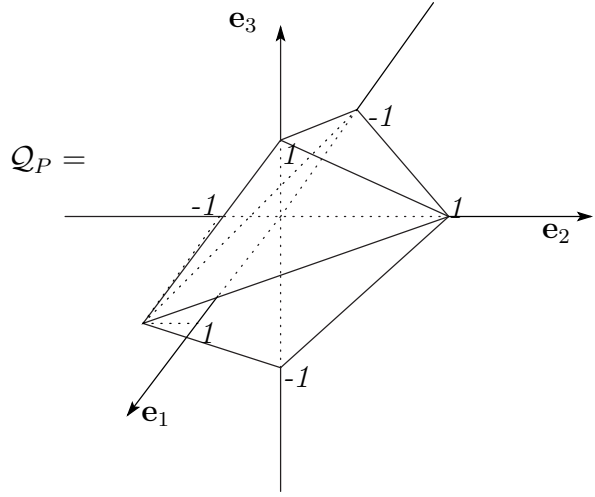
\begin{picture}( 30.4800, 25.6000)(  7.2000,-26.6000)
%
\special{pn 8}%
\special{pa 2144 820}%
\special{pa 1424 1780}%
\special{pa 3024 1220}%
\special{pa 2544 660}%
\special{pa 2144 820}%
\special{fp}%
%
\special{pn 8}%
\special{pa 2144 820}%
\special{pa 3024 1220}%
\special{fp}%
%
\special{pn 8}%
\special{pa 2144 820}%
\special{pa 2144 1516}%
\special{dt 0.045}%
%
\special{pn 8}%
\special{pa 3024 1220}%
\special{pa 1848 1220}%
\special{dt 0.045}%
%
\special{pn 8}%
\special{pa 2544 660}%
\special{pa 2144 1220}%
\special{dt 0.045}%
%
\special{pn 8}%
\special{pa 1424 1780}%
\special{pa 2144 2012}%
\special{fp}%
%
\special{pn 8}%
\special{pa 3024 1220}%
\special{pa 2144 2012}%
\special{fp}%
%
\special{pn 8}%
\special{pa 2144 1532}%
\special{pa 2144 2012}%
\special{dt 0.045}%
%
\special{pn 8}%
\special{pa 2144 1220}%
\special{pa 1824 1628}%
\special{dt 0.045}%
%
\special{pn 8}%
\special{pa 2144 2012}%
\special{pa 2144 2660}%
\special{fp}%
%
\special{pn 8}%
\special{pa 1840 1220}%
\special{pa 1016 1220}%
\special{fp}%
%
\special{pn 8}%
\special{pa 2144 820}%
\special{pa 2144 228}%
\special{fp}%
\special{sh 1}%
\special{pa 2144 228}%
\special{pa 2124 296}%
\special{pa 2144 282}%
\special{pa 2164 296}%
\special{pa 2144 228}%
\special{fp}%
%
\special{pn 8}%
\special{pa 3024 1220}%
\special{pa 3768 1220}%
\special{fp}%
\special{sh 1}%
\special{pa 3768 1220}%
\special{pa 3702 1200}%
\special{pa 3716 1220}%
\special{pa 3702 1240}%
\special{pa 3768 1220}%
\special{fp}%
%
\special{pn 8}%
\special{pa 1424 1780}%
\special{pa 2544 660}%
\special{dt 0.045}%
%
\special{pn 8}%
\special{pa 2544 660}%
\special{pa 2944 100}%
\special{fp}%
%
\special{pn 8}%
\special{pa 1816 1636}%
\special{pa 1344 2260}%
\special{fp}%
\special{sh 1}%
\special{pa 1344 2260}%
\special{pa 1400 2220}%
\special{pa 1376 2218}%
\special{pa 1368 2196}%
\special{pa 1344 2260}%
\special{fp}%
\put(18.7200,-1.5600){\makebox(0,0)[lt]{${\bf e}_3$}}%
\put(35.1200,-12.6000){\makebox(0,0)[lt]{${\bf e}_2$}}%
\put(14.1600,-22.2000){\makebox(0,0)[lt]{${\bf e}_1$}}%
\put(25.9200,-7.4800){\makebox(0,0)[lb]{-1}}%
\put(30.1600,-12.2800){\makebox(0,0)[lb]{1}}%
\put(21.4400,-8.3600){\makebox(0,0)[lt]{1}}%
\put(21.5200,-20.0400){\makebox(0,0)[lt]{-1}}%
%
\special{pn 8}%
\special{pa 1424 1780}%
\special{pa 1824 1220}%
\special{dt 0.045}%
\put(18.0800,-12.1200){\makebox(0,0)[rb]{-1}}%
%
\special{pn 8}%
\special{pa 1424 1780}%
\special{pa 1712 1780}%
\special{dt 0.045}%
\put(17.1200,-17.7200){\makebox(0,0)[lt]{1}}%
\put(7.2000,-8.5200){\makebox(0,0)[lt]{$\Xc_P=$}}%
\end{picture}%
\end{center}
\end{Example}

\bigskip

Let $P$ be a finite partially ordered set.
A subset $Q$ of $P$ is called a {\em chain}
of $P$ if $Q$ is a totally ordered subset of $P$.
The {\em length} of a chain $Q$ is $\ell(Q) = \sharp(Q) - 1$.
A chain $Q$ of $P$ is {\em saturated}
if $x, y \in Q$ with $x < y$, 
then there is no $z \in P$ with $x < z < y$.
A {\em maximal} chain of ${\hat P}$ is a saturated
chain $Q$ of ${\hat P}$ with 
$\{ {\hat 0}, {\hat 1} \} \subset Q$.

\begin{Lemma}
\label{Fano}
The convex polytope $\Xc_P$ is a Fano polytope.
\end{Lemma}

\begin{proof}
Let $e = \{ y_i, y_j \}$ be an edge of ${\hat P}$
with $y_i < y_j$.  
Let $c_e$ denote the number of maximal chains $Q$ 
of ${\hat P}$ with $\{ y_i, y_j \} \subset Q$.
If $ \{ y_{i_1}, y_{i_2}, \ldots, y_{i_m} \}$
is a maximal chain of ${\hat P}$ with
$y_0 = y_{i_1} < y_{i_2} < \ldots < y_{i_m} = y_{d+1}$, 
then
\[
\sum_{j=1}^{m-1}\rho(\{y_{i_j}, y_{i_{j+1}}\}) = (0, \ldots 0).
\]
Hence
\[
\sum_{e} c_e \rho(e) = (0, \ldots 0),
\]
where $e$ ranges all edges of ${\hat P}$.
Thus the origin of $\RR^d$ belongs to the interior of 
$\Xc_P$.  Since $\Xc_P$ is a convex polytope 
which is contained in the convex hull of the finite set  
$
\{ \sum_{i=1}^d \varepsilon_i {\bf e}_i : 
\varepsilon_i \in \{ 0, 1, -1 \} \}
$
in $\RR^d$, it follows that the origin of $\RR^d$
is the unique integer point belonging to the interior
of $\Xc_P$.  Thus $\Xc_P$ is a Fano polytope, as desired.
\end{proof}

\begin{Lemma}
\label{terminal}
The Fano polytope $\Xc_P$ is terminal.
\end{Lemma}

\begin{proof}
%

Suppose that $\Xc_P$ contains an integer point 
$\a = (\a_1,\ldots,\a_d) \in \ZZ^d$ with 
$\a \not= (0,\ldots,0).$ 
Then, obviously, $\a_1,\ldots,\a_d \in \{ -1,0,1 \}$. 
Let, say, $\a_1 = 1$. 
Let $e_1,\ldots,e_n$ be all edges of $\hat{P}$ 
and $e_{i_1},\ldots,e_{i_m}$ the edges with 
$y_1 \in e_{i_j}$ for $j = 1,\ldots,m$. 
If we set $e_{i_j} = \{ y_{i_j},y_{i_{j'}} \}$ with $ y_{i_j} < y_{i_{j'}}$, 
since $\a$ belongs to the convex hull of 
$\{ \rho(e_1),\ldots,\rho(e_n) \}$, 
then one has 
$$
\sum_{j=1}^mr_{i_j}q_{i_j} = \a_1 = 1, 
$$
where $0 \leq r_{i_1},\ldots,r_{i_m} \leq 1$ and 
$q_{i_j} = 1$ (resp. $q_{i_j} = -1$) 
if $y_1 < y_{i_{j'}}$ (resp. $y_{i_j} < y_1$). 
By removing all $r_{i_j}$ with $r_{i_j}=0$, 
we may assume that 
$$
\sum_{j=1}^{m'}r_{i_j}q_{i_j} = 1, 
$$
where $0 < r_{i_1},\ldots,r_{i_{m'}} \leq 1$. 
Since $\sum_{j=1}^{m'}r_{i_j} \leq 1$, 
there is no $j$ with $q_{i_j} = -1$. 
Hence $\sum_{j=1}^{m'}r_{i_j} = 1$. 
If $m' > 1$, then $0 < r_{i_1},\ldots,r_{i_{m'}} <1$. 
Thus $\sum_{j=1}^{m'}r_{i_j}\rho(e_{i_j})=\a \not\in \ZZ^d$. 
Thus $m' = 1$. In other words, 
if $\Xc_P$ contains an integer point $\a \not= (0,\ldots,0)$, 
then $\a$ must be one of $\rho(e_1),\ldots,\rho(e_n)$ 
and $\rho(e_1),\ldots,\rho(e_n)$ are precisely the vertices of $\Xc_P$. 
\end{proof}

\begin{Lemma}
\label{Gorenstein}
The Fano polytope $\Xc_P$ is Gorenstein.
\end{Lemma}

\begin{proof}
Via the theory of totally unimodular matrices (\cite[Chapter 9]{Sch}),
it follows that the equation of each supporting 
hyperplane of $\Xc_P$ is of the form
$a_1x_1 + \cdots + a_dx_d = 1$
with each $a_i \in \ZZ$.
In other words, the dual polytope of $\Xc_P$ is integral.
Hence $\Xc_P$ is Gorenstein, as required.
\, \, \, \, \, \, \, \, \, \, 
\, \, \, \, \, \, \, \, \, \, 
\, \, \, \, \, \, \, \, \, \, 
\, \, \, \, \, \, \, \, \, \, 
\, \, \, \, \, \, \, 
\end{proof}

\begin{Remark}{\em 
There is a well-known integral convex polytope 
arising from a finite partially ordered set $P$, 
which is called an {\em order polytope} $\Oc_P$. 
(See \cite[Chapter 4]{StanleyEC} and \cite{StanleyTPP}.) 
The authors propose to consult 
\cite{HibiDis}, \cite{Post}, \cite{Thomas} and \cite{Wagner} 
on the related work on order polytopes. 
One can verify immediately that 
the primitive outer normals of each facet of $\Oc_P$ 
one-to-one corresponds to each vertex of $\Xc_P$. 
Now $\Oc_P$ is Gorenstein if and only if $P$ is pure, 
i.e., all maximal chains of $\hat{P}$ have the same length. 
When $P$ is pure, let $l$ denote 
the length of each maximal chain of $\hat{P}$. 
Then the dilated polytope $l\Oc_P$ contains a unique integer point 
$\a \in \ZZ^d$, where $d$ is the cardinality of $P$, 
belonging to the interior of $l\Oc_P$. 
Then the dual polytope of the Gorenstein Fano polytope $l\Oc_P - \a$ 
coincides with $\Xc_P$. 
Thus, when $P$ is pure, we can associate $\Xc_P$ 
with the dual polytope of an order polytope $\Oc_P$. 
}\end{Remark}

\section{When is $\Xc_P$ $\QQ$-factorial?}
Let $P = \{ y_1, \ldots, y_d \}$
be a finite partially ordered set and
${\hat P} = P \cup \{ y_0, y_{d+1} \}$,
where $y_0 = {\hat 0}$ and $y_{d+1} = {\hat 1}$.
A sequence $\Gamma = (y_{i_1}, y_{i_2}, \ldots, y_{i_m})$ 
is called a {\em path} in ${\hat P}$ 
if $\Gamma$ is a path in the Hasse diagram of ${\hat P}$. 
In other words, 
$\Gamma = (y_{i_1}, y_{i_2}, \ldots, y_{i_m})$ 
is a path in ${\hat P}$ if 
$y_{i_j} \neq y_{i_k}$ for all $1 \leq j < k \leq m$ 
and if $\{ y_{i_j}, y_{i_{j+1}} \}$ is an edge of 
${\hat P}$ for all $1 \leq j \leq m-1$. 
In particular, if $\{ y_{i_1},y_{i_m} \}$ is 
also an edge of $\hat{P}$, 
then $\Gamma$ is called a {\em cycle}. 
The {\em length} of a path 
$\Gamma = (y_{i_1}, y_{i_2}, \ldots, y_{i_m})$ 
is $\ell(\Gamma) = m-1$ or $\ell(\Gamma)=m$ 
if $\Gamma$ is a cycle. 

A path $\Gamma = (y_{i_1}, y_{i_2}, \ldots, y_{i_{m+1}})$ 
is called {\em special} if 
\[
\sharp \{ \, j : y_{i_j} < y_{i_{j+1}}, 1 \leq j \leq m-1 \, \} 
= 
\sharp \{ \, k : y_{i_k} > y_{i_{k+1}}, 1 \leq k \leq m-1 \, \}. 
\]

Given a special path $\Gamma = (y_{i_1}, y_{i_2}, \ldots, y_{i_m})$, 
there exists a unique function 
\begin{eqnarray*}
\mu_{\Gamma} : \{ y_{i_1}, y_{i_2}, \ldots, y_{i_m} \} 
\to \{ 0, 1, 2, \ldots \} 
\end{eqnarray*}
such that
\begin{itemize}
\item
$\mu_{\Gamma}(y_{i_{j+1}}) = \mu_{\Gamma}(y_{i_j}) + 1$
(resp. $\mu_{\Gamma}(y_{i_{j}}) = \mu_{\Gamma}(y_{i_{j+1}}) + 1$)
if $y_{i_j} < y_{i_{j+1}}$
(resp. $y_{i_j} > y_{i_{j+1}}$);
\item
$\min \{ \mu_{\Gamma}(y_{i_1}), \mu_{\Gamma}(y_{i_2}), \ldots, \mu_{\Gamma}(y_{i_m}) \}
= 0$.
\end{itemize}
In particular, $\Gamma$ is special if and only if 
$\mu_{\Gamma}(y_{i_1}) = \mu_{\Gamma}(y_{i_m})$. 

Similary, a special cycle is defined and 
given a special cycle $C$, 
there exists a unique function $\mu_C$ 
which is defined the same way as above. 

\begin{Example}
{\em
Among the two paths and three cycles drawn below, 
each of one path and two cycles 
on the left-hand side is special; none of one path and one cycle 
on the right-hand side is special.

\medskip

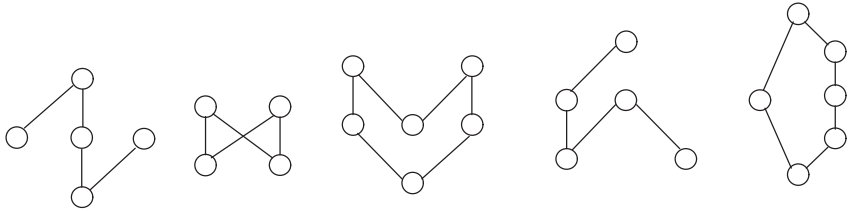
\begin{center}
\unitlength 0.1in
\begin{picture}( 43.9400, 10.7000)(  4.8600,-15.7000)
%
\special{pn 8}%
\special{ar 884 1516 56 54  0.0000000 6.2831853}%
%
\special{pn 8}%
\special{ar 542 1204 56 56  0.0000000 6.2831853}%
%
\special{pn 8}%
\special{ar 882 1204 56 56  0.0000000 6.2831853}%
%
\special{pn 8}%
\special{ar 886 896 56 54  0.0000000 6.2831853}%
%
\special{pn 8}%
\special{ar 2614 1442 56 56  0.0000000 6.2831853}%
%
\special{pn 8}%
\special{ar 2614 1138 56 54  0.0000000 6.2831853}%
%
\special{pn 8}%
\special{ar 2302 1136 56 56  0.0000000 6.2831853}%
%
\special{pn 8}%
\special{ar 2302 830 56 54  0.0000000 6.2831853}%
%
\special{pn 8}%
\special{ar 2926 1136 56 56  0.0000000 6.2831853}%
%
\special{pn 8}%
\special{ar 2926 830 56 54  0.0000000 6.2831853}%
%
\special{pn 8}%
\special{pa 2302 1084}%
\special{pa 2302 884}%
\special{fp}%
%
\special{pn 8}%
\special{pa 2662 1420}%
\special{pa 2910 1198}%
\special{fp}%
%
\special{pn 8}%
\special{ar 1530 1042 56 56  0.0000000 6.2831853}%
%
\special{pn 8}%
\special{ar 1920 1042 56 56  0.0000000 6.2831853}%
%
\special{pn 8}%
\special{ar 1920 1348 56 54  0.0000000 6.2831853}%
%
\special{pn 8}%
\special{ar 1530 1348 56 56  0.0000000 6.2831853}%
%
\special{pn 8}%
\special{pa 1570 1080}%
\special{pa 1880 1310}%
\special{fp}%
%
\special{pn 8}%
\special{pa 1920 1096}%
\special{pa 1920 1294}%
\special{fp}%
%
\special{pn 8}%
\special{ar 4826 754 56 56  0.0000000 6.2831853}%
%
\special{pn 8}%
\special{ar 4432 1008 56 54  0.0000000 6.2831853}%
%
\special{pn 8}%
\special{ar 4826 1208 56 54  0.0000000 6.2831853}%
%
\special{pn 8}%
\special{ar 4632 556 56 56  0.0000000 6.2831853}%
%
\special{pn 8}%
\special{ar 4826 984 56 56  0.0000000 6.2831853}%
%
\special{pn 8}%
\special{pa 4826 808}%
\special{pa 4826 932}%
\special{fp}%
%
\special{pn 8}%
\special{pa 4826 1038}%
\special{pa 4826 1160}%
\special{fp}%
%
\special{pn 8}%
\special{pa 1530 1094}%
\special{pa 1530 1294}%
\special{fp}%
%
\special{pn 8}%
\special{pa 1560 1310}%
\special{pa 1888 1088}%
\special{fp}%
%
\special{pn 8}%
\special{pa 2924 1084}%
\special{pa 2924 884}%
\special{fp}%
%
\special{pn 8}%
\special{pa 2326 1182}%
\special{pa 2558 1420}%
\special{fp}%
%
\special{pn 8}%
\special{pa 2332 876}%
\special{pa 2558 1114}%
\special{fp}%
%
\special{pn 8}%
\special{pa 2668 1114}%
\special{pa 2894 884}%
\special{fp}%
%
\special{pn 8}%
\special{ar 3732 702 56 56  0.0000000 6.2831853}%
%
\special{pn 8}%
\special{ar 3730 1008 56 54  0.0000000 6.2831853}%
%
\special{pn 8}%
\special{ar 3420 1008 56 56  0.0000000 6.2831853}%
%
\special{pn 8}%
\special{ar 3420 1316 56 54  0.0000000 6.2831853}%
%
\special{pn 8}%
\special{ar 4044 1316 56 54  0.0000000 6.2831853}%
%
\special{pn 8}%
\special{pa 3420 1062}%
\special{pa 3420 1262}%
\special{fp}%
%
\special{pn 8}%
\special{pa 3442 962}%
\special{pa 3676 726}%
\special{fp}%
%
\special{pn 8}%
\special{pa 3450 1268}%
\special{pa 3676 1032}%
\special{fp}%
%
\special{pn 8}%
\special{pa 3786 1032}%
\special{pa 4012 1262}%
\special{fp}%
%
\special{pn 8}%
\special{ar 4632 1398 56 56  0.0000000 6.2831853}%
%
\special{pn 8}%
\special{pa 4788 718}%
\special{pa 4666 598}%
\special{fp}%
%
\special{pn 8}%
\special{pa 4596 1364}%
\special{pa 4458 1050}%
\special{fp}%
%
\special{pn 8}%
\special{pa 4450 956}%
\special{pa 4606 598}%
\special{fp}%
%
\special{pn 8}%
\special{pa 4674 1364}%
\special{pa 4788 1254}%
\special{fp}%
%
\special{pn 8}%
\special{ar 1208 1210 56 54  0.0000000 6.2831853}%
%
\special{pn 8}%
\special{pa 582 1154}%
\special{pa 836 932}%
\special{fp}%
%
\special{pn 8}%
\special{pa 888 950}%
\special{pa 888 1148}%
\special{fp}%
%
\special{pn 8}%
\special{pa 882 1266}%
\special{pa 882 1464}%
\special{fp}%
%
\special{pn 8}%
\special{pa 928 1476}%
\special{pa 1162 1260}%
\special{fp}%
\end{picture}%
\end{center}
}
\end{Example}



We say that a path 
$\Gamma = (y_{i_1}, y_{i_2}, \ldots, y_{i_{m+1}})$
or a cycle 
$C = (y_{i_1}, y_{i_2}, \ldots, y_{i_m})$
of ${\hat P}$ {\em belongs to a facet of $\Xc_P$}
if there is a facet $\Fc$ of $\Xc_P$ with
$\rho(\{ y_{i_j}, y_{i_{j+1}} \}) \in \Fc$
for all $1 \leq j \leq m$, 
where $y_{i_{m+1}}=y_{i_1}$. 

We say that a cycle $C = (y_{i_1}, y_{i_2}, \ldots, y_{i_m})$
is {\em very special} if $C$ is special and if 
$\{ y_0,y_{d+1} \} \not\subset \{ y_{i_1}, y_{i_2}, \ldots, y_{i_m}\}$. 

\begin{Lemma}
\label{belongtoafacet}
{\em (a)}
Let $C=(y_{i_1},y_{i_2},\ldots,y_{i_m})$ be 
a cycle in ${\hat P}$. 
If $C$ belongs to a facet of $\Xc_P$, 
then $C$ is a special cycle. 
In particular, $C$ is a very special cycle or 
$C$ contains a special path 
$(y_{i_1},y_{i_2},\ldots,y_{i_{r+1}})$ with 
$y_{i_1} = y_0$ and $y_{i_{r+1}} = y_{d+1}$. 

{\em (b)}
Let $\Gamma = (y_{i_1},y_{i_2},\ldots,y_{i_m})$ with 
$y_{i_1} = y_0$ and $y_{i_m} = y_{d+1}$ 
be a path in ${\hat P}$. 
If $\Gamma$ belongs to a facet of $\Xc_P$, 
then $\Gamma$ is a special path. 
\end{Lemma}

\begin{proof}
(a) Let $a_1x_1+\cdots+a_dx_d=1$ with each $a_i \in \QQ$ 
denote the equation of the supporting hyperplane of 
$\Xc_P$ which defines the facet. 
Since $\{ y_{i_j}, y_{i_{j+1}} \}$
are edges of ${\hat P}$ for $1 \leq j \leq m$, 
where $y_{i_{m+1}} = y_{i_1}$, 
it follows that
$a_{i_j} - a_{i_{j+1}} = q_j$, 
where $q_j = 1$ if $y_{i_j} < y_{i_{j+1}}$ and 
$q_j = -1$ if $y_{i_j} > y_{i_{j+1}}$.
Now, 
\[
\sum_{i=1}^{m} q_j
= \sum_{i=1}^{m} (a_{i_j} - a_{i_{j+1}})
= 0.
\]
Hence $C$ must be special. 

Suppose that 
$\{ y_0,y_{d+1} \} \subset \{ y_{i_1}, y_{i_2}, \ldots, y_{i_m}\}$. 
Let $y_{i_1} = y_0$ and $y_{i_{r+1}} = y_{d+1}$. 
Since $\{ y_{i_j},y_{i_{j+1}} \}$, $1 \leq j \leq r$, 
are edges of $\hat{P}$, 
one has 
$-a_{i_2}=1$, $a_{i_r}=1$ and 
$a_{i_j} - a_{i_{j+1}} = q_j$ for $j=2,3,\ldots,r-1.$ 
On the one hand, one has 
$$
-a_{i_2}+\sum_{j=2}^{r-1}(a_{i_j} - a_{i_{j+1}}) + a_{i_r} = 0. 
$$
On the other hand, one has 
\begin{eqnarray*}
&&-a_{i_2}+\sum_{j=2}^{r-1}(a_{i_j} - a_{i_{j+1}}) + a_{i_r} 
= 1+\sum_{j=2}^{r-1}q_j+1 \\
&&\quad\quad= -\mu_C(y_0) + \mu_C(y_{i_2})+
\sum_{j=2}^{r-1}(\mu_C(y_{i_{j+1}}) - \mu_C(y_{i_j})) 
+\mu_C(y_{d+1}) - \mu_C(y_{i_r}) \\
&&\quad\quad= \mu_C(y_{d+1}) - \mu_C(y_0). 
\end{eqnarray*}
It then follows that one must be $\mu_C(y_0) = \mu_C(y_{d+1})$. 
Let $\Gamma = (y_{i_1}, y_{i_2},\ldots,y_{i_{r+1}}).$ 
Then it is clear that $\mu_{\Gamma}(y_0) = \mu_{\Gamma}(y_{d+1}).$ 
Thus $\Gamma$ is a special path. 
Hence $C$ contains a special path $\Gamma$. 

(b) A proof can be given by the similar way of a proof of (a). 
\end{proof}


Let $P$ be a finite partially ordered set and
$y, z \in \hat{P}$ with $y < z$.  
The {\em distance} of $y$ and $z$ in $\hat{P}$ is 
the smallest integer $s$ for which
there is a saturated chain 
$Q = \{ z_0, z_1, \ldots, z_s \}$ with
\[
y = z_0 < z_1 < \cdots < z_s = z.
\] 
Let $\dist_{\hat{P}}(y,z)$ denote the distance of $y$ and $z$
in $\hat{P}$. 


\begin{Theorem}
Let $P = \{ y_1, \ldots, y_d \}$ be a 
finite partially ordered set 
and $\hat{P} = P \union \{ y_0, y_{d+1} \}$, where
$y_0 = {\hat 0}$ and $y_{d+1} = {\hat 1}$.
Then the following conditions are equivalent:
\begin{enumerate}
\item[(i)]
$\Xc_P$ is $\QQ$-factorial;
\item[(ii)]
$\Xc_P$ is smooth;
\item[(iii)]
${\hat P}$ possesses no very special cycle 
$C = (y_{i_1}, \ldots, y_{i_m})$ 
such that
\begin{eqnarray}
\label{aaaaa}
\mu_C(y_{i_a}) - \mu_C(y_{i_b})
\leq
\dist_{\hat P}(y_{i_b}, y_{i_a}) 
\end{eqnarray}
for all $1 \leq a, \, b \leq m$ with $y_{i_b} < y_{i_a}$, and
\begin{eqnarray}
\label{bbbbb}
\mu_C(y_{i_a}) - \mu_C(y_{i_b})
\leq 
\dist_{\hat P}(y_0, y_{i_a}) + \dist_{\hat P}(y_{i_b}, y_{d+1})
\end{eqnarray}
for all $1 \leq a, \, b \leq m$, 
and no special path 
$\Gamma = (y_{i_1},\ldots,y_{i_m})$ with 
$y_{i_1} = y_0$ and $y_{i_m} = y_{d+1}$ 
such that 
\begin{eqnarray}
\label{eeeee}
\mu_{\Gamma}(y_{i_a}) - \mu_{\Gamma}(y_{i_b})
\leq
\dist_{\hat P}(y_{i_b}, y_{i_a}) 
\end{eqnarray}
for all $1 \leq a, \, b \leq m$ with $y_{i_b} < y_{i_a}$. 
\end{enumerate}
\end{Theorem}

\begin{proof}
{\bf ((i) $\Rightarrow$ (iii))}
If $C = (y_{i_1}, \ldots, y_{i_m})$
is a cycle in ${\hat P}$ with $y_{i_{m+1}} = y_1$, 
then 
\[
\sum_{j=1}^mq_j \, \rho(\{y_{i_j}, y_{i_{j+1}}\}) = (0,\ldots,0),
\]
where $q_j=1$ if $y_{i_j} < y_{i_{j+1}}$ and 
$q_j=-1$ if $y_{i_j} > y_{i_{j+1}}$. 
Thus in particular 
$\rho(\{y_{i_j}, y_{i_{j+1}}\})$,
$1 \leq j \leq m$, 
cannot be affinely independent if $C$ is special.

Now, suppose that ${\hat P}$ possesses a very special cycle 
$C = (y_{i_1}, \ldots, y_{i_m})$
which satisfies the inequalities 
(\ref{aaaaa}) and (\ref{bbbbb}).
Our work is to show that $\Xc_P$ is not simplicial.
Let $v_j = \rho(\{y_{i_j}, y_{i_{j+1}}\})$,
$1 \leq j \leq m$, where $y_{i_{m+1}} = y_{i_1}$.
Since $v_1, \ldots, v_m$ cannot be affinely independent,
to show that $\Xc_P$ is not simplicial,
what we must prove is the existence of 
a face of $\Xc_P$ which contains 
the vertices $v_1, \ldots, v_m$. 

Let $a_1, \ldots, a_d$ be integers. 
Write $\Hc \subset \RR^d$ for the hyperplane defined 
the equation $a_1x_1 + \cdots + a_dx_d = 1$
and $\Hc^{(+)} \subset \RR^d$
for the closed half-space defined by the inequality
$a_1x_1 + \cdots + a_dx_d \leq 1$.
We will determine $a_1, \ldots, a_d$ such that
$\Hc$ is a supporting hyperplane
of a face $\Fc$ of $\Xc_P$ with
$\{ v_1, \ldots, v_m \} \subset \Fc$
and with $\Xc_P \subset \Hc^{(+)}$.

{\em First Step.}
It follows from $(\ref{bbbbb})$ that
\begin{eqnarray}
\label{ccccc}
\max_{1 \leq a \leq m} (\mu_C(y_{i_a}) 
- \dist_{\hat P}(y_0, y_{i_a}))
\leq 
\min_{1 \leq b \leq m} (\mu_C(y_{i_b}) 
+ \dist_{\hat P}(y_{i_b}, y_{d+1})). 
\end{eqnarray}
By using $(\ref{aaaaa})$, 
if $y_0 \in \{ y_{i_1}, \ldots, y_{i_m} \}$, then
the left-hand side of $(\ref{ccccc})$
is equal to $\mu_C(y_0)$.  Similarly, if
$y_{d+1} \in \{ y_{i_1}, \ldots, y_{i_m} \}$, then
the right-hand side of $(\ref{ccccc})$
is equal to $\mu_C(y_{d+1})$.

Now, fix an arbitrary integer $a$ with
\[
\max_{1 \leq a \leq m} (\mu_C(y_{i_a}) 
- \dist_{\hat P}(y_0, y_{i_a}))
\leq a \leq 
\min_{1 \leq b \leq m} (\mu_C(y_{i_b}) 
+ \dist_{\hat P}(y_{i_b}, y_{d+1})). 
\]
However, exceptionally,  
if $y_0 \in \{ y_{i_1}, \ldots, y_{i_m} \}$,
then $a = \mu_C(y_0)$. 
If $y_{d+1} \in \{ y_{i_1}, \ldots, y_{i_m} \}$, 
then 
$a = \mu_C(y_{d+1})$.
Let $a_{i_j} = a - \mu_C({y_{i_j}})$
for $1 \leq j \leq m$.
Then one has
\begin{eqnarray}
\label{qqqqq}
- a_{i_j} \leq \dist_{\hat P}(y_0, y_{i_j}), \, \, \, \, \, 
a_{i_j} \leq \dist_{\hat P}(y_{i_j}, y_{d+1}).
\end{eqnarray}
Moreover, it follows easily that 
each $v_j$ lies on the hyperplane of $\RR^d$ 
defined by the equation
\[
\sum_{{i_j} \not\in \{ 0, \, d+1 \}}
a_{i_j} x_{i_j} = 1.
\]

{\em Second Step.}
Let $A = {\hat P} 
\setminus (\{ y_0, y_{d+1} \} \cup \{ y_{i_1}, \ldots, y_{i_m} \})$
and $y_i \in A$.
\begin{itemize}
\item
Suppose that there is $y_{i_j}$ with $y_{i_j} < y_i$ and that
there is no $y_{i_k}$ with $y_{i_k} > y_i$.
Then we define $a_i$ by setting
\[
a_i = 
\max (\{ a_{i_j} - \dist_{\hat P}(y_{i_j}, y_i) : y_{i_j} < y_i \}
\cup \{ 0 \}).
\]
\item
Suppose that there is no $y_{i_j}$ with $y_{i_j} < y_i$ and that
there is $y_{i_k}$ with $y_{i_k} > y_i$.
Then we define $a_i$ by setting
\[
a_i = 
\min (\{ a_{i_k} + \dist_{\hat P}(y_{i}, y_{i_k}) : y_{i} < y_{i_k} \}
\cup \{ 0 \}).
\]
\item
Suppose that there is $y_{i_j}$ with $y_{i_j} < y_i$ and that
there is $y_{i_k}$ with $y_{i_k} > y_i$.
Then either 
\[
b_i = 
\max (\{ a_{i_j} - \dist_{\hat P}(y_{i_j}, y_i) : y_{i_j} < y_i \}
\cup \{ 0 \})
\]
or
\[
c_i = 
\min (\{ a_{i_k} + \dist_{\hat P}(y_{i}, y_{i_k}) : y_{i} < y_{i_k} \}
\cup \{ 0 \})
\]
must be zero.  In fact, if $b_i \neq 0$ and $c_i \neq 0$,
then there are $j$ and $k$ with
$a_{i_j} > \dist_{\hat P}(y_{i_j}, y_i)$ and 
$- a_{i_k} > \dist_{\hat P}(y_{i}, y_{i_k})$.
Since $\mu_C(y_{i_k}) - \mu_C(y_{i_j})
= a_{i_j} - a_{i_k}$ and since 
$\dist_{\hat P}(y_{i_j}, y_i) + \dist_{\hat P}(y_{i}, y_{i_k})
\geq \dist_{\hat P}(y_{i_j}, y_{i_k})$,
it follows that
\begin{eqnarray*}
\mu_C(y_{i_k}) - \mu_C(y_{i_j})
> \dist_{\hat P}(y_{i_j}, y_{i_k}).
\end{eqnarray*}
This contradicts $(\ref{aaaaa})$.  Hence either $b_i = 0$
or $c_i = 0$.  If $b_i \neq 0$, then we set $a_i = b_i$.
If $c_i \neq 0$, then we set $a_i = c_i$.
If $b_i = c_i = 0$, then we set $a_i = 0$.
\item
Suppose that there is no $y_{i_j}$ with $y_{i_j} < y_i$ and that
there is no $y_{i_k}$ with $y_{i_k} > y_i$.
Then we set $a_i = 0$. 
\end{itemize}

{\em Third Step.}
Finally, we finish determining the integers $a_1, \ldots, a_d$.
Let $\Hc \subset \RR^d$ denote the hyperplane 
defined by the equation
$a_1x_1 + \ldots + a_dx_d = 1$ and
$\Hc^{(+)} \subset \RR^d$ the closed half-space defined by
the inequality $a_1x_1 + \ldots + a_dx_d \leq 1$.
Since each $v_j$ lies on the hyperplane $\Hc$, 
in order for $\Fc = \Hc \cap \Xc_P$ to be a face of $\Xc_P$,
it is required to show $\Xc_P \subset \Hc^{(+)}$. 
Let $\{ y_i, y_j \}$ with $y_i < y_j$ be an edge of ${\hat P}$.
\begin{itemize}
\item
Let $y_i \in \{ y_{i_1}, \ldots, y_{i_m} \}$ with
$y_j \not\in \{ y_{i_1}, \ldots, y_{i_m} \}$.
If $y_j \neq y_{d+1}$, then
\[a_j \geq \max \{ a_i - 1, 0 \},\]  
where $a_0 = 0$.  Thus 
$a_i - a_j \leq 1$.
If $y_j = y_{d+1}$, then by using $(\ref{qqqqq})$ 
one has $a_i \leq 1$, as desired. 
\item
Let $y_j \in \{ y_{i_1}, \ldots, y_{i_m} \}$ with
$y_i \not\in \{ y_{i_1}, \ldots, y_{i_m} \}$.
If $y_i \neq y_0$, then 
\[a_i \leq \min \{ a_j + 1, 0 \},\]  
where $a_{d+1} = 0$.  Thus
$a_i - a_j \leq 1$. 
If $y_i = y_{0}$, then by using $(\ref{qqqqq})$ 
one has $- a_j \leq 1$, as desired. 
\end{itemize}

Let $A' = {\hat P} 
\setminus \{ y_{i_1}, \ldots, y_{i_m} \}$.
Write $B$ for the subset of $A'$ consisting of those
$y_i \in A'$ such that there is $j$ with $y_{i_j} < y_i$.
Write $C$ for the subset of $A'$ consisting of those
$y_i \in A'$ such that there is $k$ with $y_i < y_{i_k}$.
Again, 
let $e = \{ y_i, y_j \}$ with $y_i < y_j$ be an edge of ${\hat P}$.
In each of the nine cases below, a routine computation
easily yields that $\rho(e) \in \Hc^{(+)}$. 
\begin{itemize}
\item
$y_i \in B \setminus C$ and $y_j \in B \setminus C$;
\item
$y_i \in C \setminus B$ and $y_j \in C \setminus B$;
\item
$y_i \in C \setminus B$ and $y_j \in B \setminus C$;
\item
$y_i \in C \setminus B$ and $y_j \in B \cap C$;
\item
$y_i \in C \setminus B$ and $y_j \not\in B \cup C$;
\item
$y_i \in B \cap C$ and $y_j \in B \cap C$; 
\item
$y_i \in B \cap C$ and $y_j \in B \setminus C$;
\item
$y_i \not\in B \cup C$ and $y_j \in B \setminus C$; 
\item
$y_i \not\in B \cup C$ and $y_j \not\in B \cup C$. 
\end{itemize}
For example, in the first case, a routine computation
is as follows.  
Let $y_j \neq y_{d+1}$.  
Let $a_i = 0$.  Then, since $a_j \geq 0$, one has $a_i - a_j \leq 1$.
Let $a_i > 0$.  Then, since $a_j \geq a_i - 1$, one has
$a_i - a_j \leq 1$.  Let $y_j = y_{d+1}$ and $a_i > 0$.
Then there is $j$ with 
$a_i = a_{i_j} - \dist_{\hat P}(y_{i_j}, y_i)$.
By using $(\ref{qqqqq})$ one has
$a_{i_j} \leq \dist_{\hat P}(y_{i_j}, y_{d+1})$.
Thus $a_i \leq  
\dist_{\hat P}(y_{i_j}, y_{d+1}) - \dist_{\hat P}(y_{i_j}, y_i)$.
Hence $a_i \leq 1$, as required.

{\em Fourth step.} 
Suppose that ${\hat P}$ possesses a special path 
$\Gamma = (y_{i_1},y_{i_2},\ldots,y_{i_m})$ 
with $y_{i_1} = y_0$ and $y_{i_m} = y_{d+1}$ 
which satisfies the inequalities (\ref{eeeee}). 
Then one has 
\[
\sum_{j=1}^{m-1}q_j \, \rho(\{y_{i_j}, y_{i_{j+1}}\}) = (0,\ldots,0),
\]
where $q_j=1$ if $y_{i_j} < y_{i_{j+1}}$ and 
$q_j=-1$ if $y_{i_j} > y_{i_{j+1}}$. 
Thus $\rho(\{y_{i_j}, y_{i_{j+1}}\})$, 
$1 \leq j \leq m-1$, 
cannot be affinely independent. 
Our work is to show that $\Xc_P$ is not simplicial. 
In this case, however, the same discussion can be given 
as the case which ${\hat P}$ possesses a very special cycle. 
(We should set $a = \mu_{\Gamma}(y_0)$ ($=\mu_{\Gamma}(y_{d+1})$).) 

\smallskip

{\bf ((iii) $\Rightarrow$ (i))}
Now, suppose that $\Xc_P$ is not $\QQ$-factorial. 
Thus $\Xc_P$ possesses a facet $\Fc$ which is not a simplex. 
Let $v_1, \ldots, v_n$ denote the vertices of $\Fc$,
where $n > d$, and
$e_j$ the edge of ${\hat P}$ with $v_j = \rho(e_j)$ 
for $1 \leq j \leq n$.
Let $a_1x_1 + \cdots + a_dx_d = 1$ denote the equation
of the supporting hyperplane $\Hc \subset \RR^d$ of $\Xc_P$
with $\Fc = \Xc_P \cap \Hc$ and 
with $\Xc_P \subset \Hc^{(+)}$, where 
$\Hc^{(+)} \subset \RR^d$
is the closed-half space
defined by the inequality  
$a_1x_1 + \cdots + a_dx_d \leq 1$.
Since $v_1, \ldots, v_n$ are not affinely independent,
there is $(r_1, \ldots, r_n) \in \ZZ^n$ with
$(r_1, \ldots, r_n) \neq (0, \ldots, 0)$ 
such that
$r_1v_1+\cdots+r_nv_n=(0,\ldots,0)$. 
By removing $r_j$ with $r_j = 0$, we may assume that
$r_1v_1+\cdots+r_{n'}v_{n'}=(0,\ldots,0)$, where
$r_j \neq 0$ for $1 \leq j \leq n'$
with $r_1 + \cdots + r_{n'} = 0$. 
Let $e_j = \{ y_{i_j}, y_{i_{j'}} \}$ with $1 \leq i_j,i_{j'} \leq d$. 
If either $y_{i_j}$ or $y_{i_{j'}}$ 
appears only in $e_j$ among
the edges $e_1, \ldots, e_{n'}$, then
$r_j = 0$.  Hence both 
$y_{i_j}$ and $y_{i_{j'}}$ must appear 
in at least two edges among 
$e_1, \ldots, e_{n'}$. 
Let $G$ denote the subgraph of the Hasse diagram 
of ${\hat P}$ with the edges
$e_1, \ldots, e_{n'}$.
Then there is no end point of $G$ in $P$. 
Thus $G$ possesses a cycle of $\hat{P}$ or 
$G$ is a path of $\hat{P}$ 
from $y_0$ to $y_{d+1}$. 
Since $v_1,\ldots,v_{n'}$ 
are contained in the facet $\Fc$,
Lemma \ref{belongtoafacet} says that
every cycle in $G$ is very special or else 
$G$ contains a special path. 

Suppose that $G$ possesses a very special cycle 
$C =(y_{i_1}, y_{i_2}, \ldots, y_{i_m})$. 
Our goal is to show that $C$ satisfies 
the inequalities $(\ref{aaaaa})$ and $(\ref{bbbbb})$.

Let  
$y_{k_0} < y_{k_1} < \cdots < y_{k_\ell}$
be a saturated chain of ${\hat P}$ 
with $\ell = \dist_{\hat P}(y_{k_0}, y_{k_\ell})$ 
such that each of  
$y_{k_0}$ and $y_{k_\ell}$ belongs to
$\{ y_{i_1}, y_{i_2}, \ldots, y_{i_m} \}$.
We claim \[
\mu_C(y_{k_\ell}) - \mu_C(y_{k_0})
\leq \dist_{\hat P}(y_{k_0}, y_{k_{\ell}}).
\]
\begin{itemize}
\item
Let $y_0 \neq y_{k_0}$ and $y_{d+1} \neq y_{k_\ell}$. 
Since ${\bf e}_{k_j} - {\bf e}_{k_{j+1}} \in \Xc_P$,
one has $a_{k_j} - a_{k_{j+1}} \leq 1$
for each $0 \leq j \leq \ell - 1$. 
Hence $a_{k_0} - a_{k_{\ell}} \leq \ell$. 
On the other hand, 
$a_{k_0} - a_{k_{\ell}} = 
\mu_C(y_{k_\ell}) - \mu_C(y_{k_0})$.
Thus $\mu_C(y_{k_\ell}) - \mu_C(y_{k_0}) 
\leq \dist_{\hat P}(y_{k_0}, y_{k_{\ell}})$.
\item
Let $y_0 = y_{k_0}$ and $y_{d+1} \neq y_{k_\ell}$. 
Since $-{\bf e}_{k_1} \in \Xc_P$,
one has $- a_{k_1} \leq 1$.
Since ${\bf e}_{k_j} - {\bf e}_{k_{j+1}} \in \Xc_P$,
one has $a_{k_j} - a_{k_{j+1}} \leq 1$
for each $1 \leq j \leq \ell - 1$. 
Hence $a_{k_1} - a_{k_{\ell}} \leq \ell - 1$. 
Thus $- a_{k_{\ell}} \leq \ell$.
On the other hand, 
$- a_{k_{\ell}} = 
\mu_C(y_{k_\ell}) - \mu_C(y_{k_0})$. 
Thus $\mu_C(y_{k_\ell}) - \mu_C(y_{k_0}) 
\leq \dist_{\hat P}(y_{k_0}, y_{k_{\ell}})$.
\item
Let $y_0 \neq y_{k_0}$ and $y_{d+1} = y_{k_\ell}$. 
Since ${\bf e}_{k_j} - {\bf e}_{k_{j+1}} \in \Xc_P$,
one has $a_{k_j} - a_{k_{j+1}} \leq 1$
for each $0 \leq j \leq \ell - 2$. 
Hence $a_{k_0} - a_{k_{\ell-1}} \leq \ell - 1$. 
Since ${\bf e}_{k_{\ell-1}} \in \Xc_P$,
one has $a_{k_{\ell-1}} \leq 1$.
Hence $a_{k_0} \leq \ell$.
On the other hand, 
$a_{k_0} = 
\mu_C(y_{k_\ell}) - \mu_C(y_{k_0})$.
Thus $\mu_C(y_{k_\ell}) - \mu_C(y_{k_0}) 
\leq \dist_{\hat P}(y_{k_0}, y_{k_{\ell}})$.
\end{itemize}

Finally, fix arbitrary $y_{i_j}$ and $y_{i_k}$ 
with 
$\mu_C(y_{i_j}) < \mu_C(y_{i_k})$.
Then 
$- a_{i_k} \leq
\dist_{\hat P}(y_0, y_{i_k})$
and
$a_{i_j} \leq
\dist_{\hat P}(y_{i_j}, y_{d+1})$.
We claim
\[
\mu_C(y_{i_k}) - \mu_C(y_{i_j})
\leq \dist_{\hat P}(y_0, y_{i_k}) + \dist_{\hat P}(y_{i_j}, y_{d+1}).
\]
If $y_{i_j} \neq y_0$ and $y_{i_k} \neq y_{d+1}$,
then  
$a_{i_j} - a_{i_k} = \mu_C(y_{i_k}) - \mu_C(y_{i_j})$. 
If $y_{i_j} = y_0$ and $y_{i_k} \neq y_{d+1}$,
then $- a_{i_k} = \mu_C(y_{i_k}) - \mu_C(y_{i_j})$. 
If $y_{i_j} \neq y_0$ and $y_{i_k} = y_{d+1}$,
then  
$a_{i_j} = \mu_C(y_{i_k}) - \mu_C(y_{i_j})$. 
Hence the required inequality follows immediately.

Suppose that $G$ contains a special path 
$\Gamma =(y_{i_1}, y_{i_2}, \ldots, y_{i_m})$ with 
$y_{i_1} = y_0$ and $y_{i_m} = y_{d+1}$. 
Our goal is to show that $C$ satisfies 
the inequalities $(\ref{eeeee})$. 
Now the same discussion can be given as above. 

\smallskip

{\bf ((i) $\Rightarrow$ (ii))}
If $P$ is a totally ordered set, then
$\Xc_P$ is a $d$-simplex with the vertices, say,
$- {\bf e}_1, {\bf e}_1 - {\bf e}_2, \ldots, 
{\bf e}_{d-1} - {\bf e}_d, {\bf e}_d$.
Thus in particular $\Xc_P$ is smooth.

Now, suppose that $P$ is not a totally ordered set.
Then ${\hat P}$ possesses a cycle.
Let $C = (y_{i_1}, \ldots, y_{i_m})$ be a cycle in ${\hat P}$.
If $C$ is not special, then Lemma \ref{belongtoafacet} (a)
says that $C$ cannot belong to a facet of $\Xc_P$.
If $C$ is special, then as was shown in the proof of
(i) $\Rightarrow$ (iii) it follows that
$\rho(\{y_{i_j}, y_{i_{j+1}}\})$,
$1 \leq j \leq m$, where $y_{i_{m+1}} = y_{i_1}$,
are not affinely independent.
Hence there is no facet $\Fc$ of $\Xc_P$ with
$\rho(\{y_{i_j}, y_{i_{j+1}}\}) \in \Fc$
for all $1 \leq j \leq m$. 

Let $\Fc$ be an arbitrary facet of $\Xc_P$ 
with $d$ vertices
$v_j = \rho(e_j)$, $1 \leq j \leq d$.
Let $G$ denote the subgraph of the Hasse diagram
of ${\hat P}$ with the edges
$e_1, \ldots, e_d$ 
and $V(G)$ the vertex set of $G$.
Since $\Fc$ is of dimension $d-1$, 
it follows that, for each $1 \leq i \leq d$, 
there is a vertex of $\Fc$ whose
$i$th coordinate is nonzero.
Hence $P \subset V(G)$.
Suppose that $P = V(G)$.
Since $G$ has $d$ edges, it follows that
$G$ possesses a cycle, a contradiction.
Hence either $y_0 \in V(G)$ or $y_{d+1} \in V(G)$.

What we must prove is that the determinant
\begin{eqnarray}
\label{determinant}
\begin{vmatrix}
v_1 \\
\vdots \\
v_d 
\end{vmatrix}
\end{eqnarray}
is equal to $\pm 1$.
Let, say, $e_1 = \{ y_1, y_{d+1} \}$.
Thus $v_1 = (1,0,\ldots,0)$.
Now, since $G$ is a forest, 
by arranging the numbering of the elements of $P$
if necessary,   
one has 
\[
\begin{vmatrix}
v_1 \\
\vdots \\
v_d 
\end{vmatrix}
=
\begin{vmatrix}
&a_{11}   &0     &\cdots  &\cdots  &0         \\
&a_{21}   &a_{22}&\ddots  &\ddots  &\vdots    \\
&\vdots   &\ddots&\ddots  &\ddots  &\vdots    \\
&\vdots   &\ddots&\ddots  &\ddots  &0         \\
&a_{d1}   &a_{d2}&\cdots  &\cdots  &a_{dd}
\end{vmatrix}, 
\]
with each $a_{ij} \in \{1,0,-1\}$. 
Since the determinant $(\ref{determinant})$ is nonzero,
it follows that
the determinant $(\ref{determinant})$ is equal to
$\pm 1$,
as desired.

\smallskip

{\bf ((ii) $\Rightarrow$ (i))}
In general, every smooth Fano polytope is 
$\QQ$-factorial.
\, \, \, \, \, \, 
\end{proof}

\begin{Corollary}
Suppose that a finite partially ordered set $P$ is pure.
Then the following conditions are equivalent:
\begin{enumerate}
\item[(i)]
$\Xc_P$ is $\QQ$-factorial;
\item[(ii)]
$\Xc_P$ is smooth;
\item[(iii)]
$P$ is a disjoint union of chains.  
\end{enumerate}
\end{Corollary}

\begin{proof}
If $P$ is pure, then
every cycle of ${\hat P}$ is special and, in addition,
satisfies the inequalities 
$(\ref{aaaaa})$ and $(\ref{bbbbb})$. 
Moreover, every path from $y_0$ to $y_{d+1}$ cannot be special. 
Hence $\Xc_P$ is $\QQ$-factorial if and only if
there is no very special cycle, i.e., 
every cycle of ${\hat P}$ possesses both ${\hat 0}$ 
and ${\hat 1}$. 
Now if there is a connected component of $P$ 
which is not a chain, 
then $P$ possesses a very special cycle. 
Thus $\Xc_P$ is $\QQ$-factorial if and only if 
$P$ does not possess a connected component 
which is not a chain. 
In other words, 
$\Xc_P$ is $\QQ$-factorial if and only if
$P$ is a disjoint union of chains, as desired. 
\, \, \, \, \, \, \, \, \, \, 
\, \, \, \, \, \, \, \, \, \, 
\, \, \, \, \, \, \, \, \, \, 
\end{proof}

\begin{Example}
{\em
Among the five finite partially ordered sets drawn below, 
each of the three finite partially ordered sets 
on the left-hand side yields a $\QQ$-factorial Fano polytope; 
none of the two finite partially ordered sets 
on the right-hand side yields a $\QQ$-factorial Fano polytope.
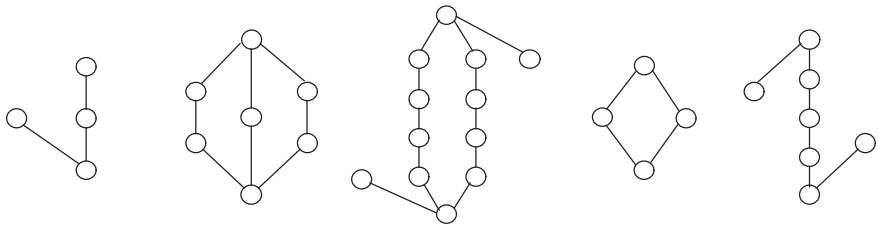
\begin{center}
\unitlength 0.1in
\begin{picture}( 45.4700, 11.3900)(  8.3300,-16.7000)
%
\special{pn 8}%
\special{ar 4170 1392 52 48  0.0000000 6.2831853}%
%
\special{pn 8}%
\special{ar 3952 1114 52 48  0.0000000 6.2831853}%
%
\special{pn 8}%
\special{pa 4126 1364}%
\special{pa 3974 1160}%
\special{fp}%
%
\special{pn 8}%
\special{ar 4390 1120 52 50  0.0000000 6.2831853}%
%
\special{pn 8}%
\special{ar 4172 844 52 48  0.0000000 6.2831853}%
%
\special{pn 8}%
\special{pa 4206 1350}%
\special{pa 4346 1154}%
\special{fp}%
%
\special{pn 8}%
\special{pa 3974 1072}%
\special{pa 4126 876}%
\special{fp}%
%
\special{pn 8}%
\special{ar 2116 708 52 50  0.0000000 6.2831853}%
%
\special{pn 8}%
\special{ar 2114 1114 54 48  0.0000000 6.2831853}%
%
\special{pn 8}%
\special{ar 1824 980 52 48  0.0000000 6.2831853}%
%
\special{pn 8}%
\special{ar 1824 1250 52 50  0.0000000 6.2831853}%
%
\special{pn 8}%
\special{ar 2408 980 52 48  0.0000000 6.2831853}%
%
\special{pn 8}%
\special{ar 2408 1250 52 50  0.0000000 6.2831853}%
%
\special{pn 8}%
\special{pa 1824 1026}%
\special{pa 1824 1202}%
\special{fp}%
%
\special{pn 8}%
\special{pa 2160 728}%
\special{pa 2392 924}%
\special{fp}%
%
\special{pn 8}%
\special{ar 886 1120 52 50  0.0000000 6.2831853}%
%
\special{pn 8}%
\special{ar 1250 1120 52 50  0.0000000 6.2831853}%
%
\special{pn 8}%
\special{ar 1250 1392 52 48  0.0000000 6.2831853}%
%
\special{pn 8}%
\special{ar 1250 850 52 48  0.0000000 6.2831853}%
%
\special{pn 8}%
\special{pa 922 1156}%
\special{pa 1212 1358}%
\special{fp}%
%
\special{pn 8}%
\special{pa 1250 898}%
\special{pa 1250 1074}%
\special{fp}%
%
\special{pn 8}%
\special{pa 1250 1168}%
\special{pa 1250 1344}%
\special{fp}%
%
\special{pn 8}%
\special{pa 4352 1080}%
\special{pa 4214 870}%
\special{fp}%
%
\special{pn 8}%
\special{pa 2406 1026}%
\special{pa 2406 1202}%
\special{fp}%
\special{pa 2114 754}%
\special{pa 2114 1066}%
\special{fp}%
%
\special{pn 8}%
\special{ar 2114 1520 54 50  0.0000000 6.2831853}%
%
\special{pn 8}%
\special{pa 1860 1282}%
\special{pa 2078 1486}%
\special{fp}%
%
\special{pn 8}%
\special{pa 2152 1486}%
\special{pa 2370 1282}%
\special{fp}%
%
\special{pn 8}%
\special{pa 2114 1160}%
\special{pa 2114 1472}%
\special{fp}%
%
\special{pn 8}%
\special{ar 5036 708 54 50  0.0000000 6.2831853}%
%
\special{pn 8}%
\special{ar 5036 1120 52 48  0.0000000 6.2831853}%
%
\special{pn 8}%
\special{ar 4746 980 52 48  0.0000000 6.2831853}%
%
\special{pn 8}%
\special{ar 5328 1250 52 50  0.0000000 6.2831853}%
%
\special{pn 8}%
\special{ar 5036 1520 52 50  0.0000000 6.2831853}%
%
\special{pn 8}%
\special{pa 5072 1486}%
\special{pa 5290 1282}%
\special{fp}%
%
\special{pn 8}%
\special{ar 5036 916 52 50  0.0000000 6.2831853}%
%
\special{pn 8}%
\special{ar 5036 1324 52 48  0.0000000 6.2831853}%
%
\special{pn 8}%
\special{pa 5036 762}%
\special{pa 5036 870}%
\special{fp}%
%
\special{pn 8}%
\special{pa 5036 964}%
\special{pa 5036 1072}%
\special{fp}%
%
\special{pn 8}%
\special{pa 5036 1168}%
\special{pa 5036 1276}%
\special{fp}%
%
\special{pn 8}%
\special{pa 5036 1372}%
\special{pa 5036 1480}%
\special{fp}%
%
\special{pn 8}%
\special{ar 2992 810 52 48  0.0000000 6.2831853}%
%
\special{pn 8}%
\special{ar 2992 1222 52 48  0.0000000 6.2831853}%
%
\special{pn 8}%
\special{ar 2992 1020 52 50  0.0000000 6.2831853}%
%
\special{pn 8}%
\special{ar 2992 1426 52 50  0.0000000 6.2831853}%
%
\special{pn 8}%
\special{pa 2992 864}%
\special{pa 2992 972}%
\special{fp}%
%
\special{pn 8}%
\special{pa 2992 1066}%
\special{pa 2992 1176}%
\special{fp}%
%
\special{pn 8}%
\special{pa 2992 1270}%
\special{pa 2992 1378}%
\special{fp}%
%
\special{pn 8}%
\special{ar 3290 810 52 48  0.0000000 6.2831853}%
%
\special{pn 8}%
\special{ar 3290 1222 52 48  0.0000000 6.2831853}%
%
\special{pn 8}%
\special{ar 3290 1020 52 50  0.0000000 6.2831853}%
%
\special{pn 8}%
\special{ar 3290 1426 52 50  0.0000000 6.2831853}%
%
\special{pn 8}%
\special{pa 3290 864}%
\special{pa 3290 972}%
\special{fp}%
%
\special{pn 8}%
\special{pa 3290 1066}%
\special{pa 3290 1176}%
\special{fp}%
%
\special{pn 8}%
\special{pa 3290 1270}%
\special{pa 3290 1378}%
\special{fp}%
%
\special{pn 8}%
\special{ar 3136 580 52 48  0.0000000 6.2831853}%
%
\special{pn 8}%
\special{ar 3136 1622 52 48  0.0000000 6.2831853}%
%
\special{pn 8}%
\special{ar 3572 810 54 48  0.0000000 6.2831853}%
%
\special{pn 8}%
\special{ar 2692 1440 52 50  0.0000000 6.2831853}%
%
\special{pn 8}%
\special{pa 3004 762}%
\special{pa 3098 606}%
\special{fp}%
%
\special{pn 8}%
\special{pa 3178 612}%
\special{pa 3274 768}%
\special{fp}%
%
\special{pn 8}%
\special{pa 3186 586}%
\special{pa 3536 774}%
\special{fp}%
%
\special{pn 8}%
\special{pa 3018 1466}%
\special{pa 3098 1602}%
\special{fp}%
%
\special{pn 8}%
\special{pa 3178 1588}%
\special{pa 3260 1458}%
\special{fp}%
%
\special{pn 8}%
\special{pa 3084 1616}%
\special{pa 2734 1458}%
\special{fp}%
%
\special{pn 8}%
\special{pa 4766 930}%
\special{pa 4992 728}%
\special{fp}%
%
\special{pn 8}%
\special{pa 1854 938}%
\special{pa 2056 728}%
\special{fp}%
\end{picture}%
\end{center}
}\end{Example}

Let $P$ and $P'$ be finite partially ordered sets. 
Then one can verify easily that 
$\Xc_P$ is isomorphic with $\Xc_{P'}$ as a convex polytope if and only if 
$P$ is isomorphic with $P'$ or with 
the dual finite partially ordered set of $P'$ 
as a finite partially ordered set. 

On the following table drawn below, 
the number of finite partially ordered sets with $d (\leq 8)$ elements, 
up to isomorphic and up to isomorphic with dual finite partially ordered sets, 
is written in the second row. Moreover, among those, 
the number of finite partially ordered sets 
constructing smooth Fano polytopes is written in the third row. 
\begin{center}
\begin{tabular}{|c|c|c|c|c|c|c|c|c|} \hline
       &$d=1$ &$d=2$ &$d=3$ &$d=4$ &$d=5$ &$d=6$ &$d=7$ &$d=8$ \\ \hline
Posets &1     &2     &4     &12    &39    &184   &1082  &8746  \\ \hline
Smooth &1     &2     &3     &6     &12    &31    &83    &266   \\ \hline
\end{tabular}
\end{center}

\end{document}